# Improvement on Brook's theorem for $3K_1$–free Graphs

Medha Dhurandhar

**Abstract:** Problem of finding an optimal upper bound for the chromatic no. of a $3K_1$-free graph is still open and pretty hard. Here we prove that for a $3K_1$-free graph G with $\Delta(G) \geq 8$, $\chi(G) \leq \max\{\Delta-1, \omega\}$. We also prove that if G is $3K_1$-free, $\omega = 4$ and $\Delta(G) \geq 7$, then $\chi(G) \leq \Delta-1$. This implies that **Borodin & Kostochka Conjecture** is true for $3K_1$-free graphs as a corollary.

**Introduction:**

In [1], [2], [3], [4] chromatic bounds for graphs are considered especially in relation with $\omega$ and $\Delta$. Gyárfás [5] and Kim [6] show that the optimal $\chi$-binding function for the class of $3K_1$-free graphs has order $\omega^2/\log(\omega)$. If we forbid additional induced subgraphs, the order of the optimal $\chi$-binding function drops below $\omega^2/\log(\omega)$. In 1941, Brooks' theorem stated that for any connected undirected graph $G$ with maximum degree $\Delta$, the chromatic number of $G$ is at most $\Delta$ unless $G$ is a complete graph or an odd cycle, in which case the chromatic number is $\Delta + 1$ [5]. In 1977, **Borodin & Kostochka** conjectured that if $\Delta(G) \geq 9$, then $\chi(G) \leq \max\{\omega, \Delta-1\}$ [6]. In 1999, Reed proved the conjecture for $\Delta \geq 10^{14}$ [7]. Also D. W. Cranston and L. Rabern [8] proved it for claw-free graphs. Here we prove that if a graph G is $3K_1$-free and $\Delta(G) \geq 8$, then $\chi(G) \leq \max\{\Delta-1, \omega\}$. For a $3K_1$-free graph with $\omega = 4$ we prove a stronger result that if $\Delta(G) \geq 7$, then $\chi(G) \leq \Delta-1$. These results prove **Borodin & Kostochka** conjecture for $3K_1$-free graphs as a corollary.

**Notation:** For a graph G, V(G), E(G), $\Delta$, $\omega$, $\chi$ denote the vertex set, edge set, maximum degree, size of a maximum clique, chromatic number of G resply. For $u \in V(G)$, $N(u) = \{v \in V(G) / uv \in E(G)\}$, and $\overline{N(u)} = N(u) \cup (u)$. If $S \subseteq V$, then $<S>$ denotes the subgraph of G induced by S. If C is some coloring of G and if a vertex u of G is colored m in C, then u is called a m-vertex. Also if P is a path in G s.t. vertices on P are alternately colored say i and j, then P is called an i-j path. All graphs considered henceforth are simple. We consider here simple and undirected graphs. For terms which are not defined herein we refer to Bondy and Murty [9].

**Main Result 1:** Let G be $3K_1$-free, if $\omega = 4$ and $\Delta \geq 7$, then $\chi \leq \Delta-1$.
Proof: Let if possible G be a smallest $3K_1$-free graph with $\Delta \geq 7$ and $\chi > \Delta-1$. Then clearly as $G \neq C_{2n+1}$ or $K_{|V(G)|}$, $\chi = \Delta > \omega$. Let $u \in V(G)$. Then $G-u \neq K_{|V(G)|-1}$ (else $\chi = \omega$). If $\Delta(G-u) \geq 7$, then by minimality $\chi(G-u) \leq \max\{\omega(G-u), \Delta(G-u)-1\}$. Clearly if $\omega(G-u) \leq \Delta(G-u)-1$, then $\chi(G-u) = \Delta(G-u)-1 \leq \Delta-1$ and otherwise $\chi(G-u) = \omega(G-u) \leq \omega < \Delta$. In any case $\chi(G-u) \leq \Delta-1$. Also if $\Delta(G-u) < 7$, then as $G-u \neq C_{2n+1}$ (else as G is $3K_1$-free, $G-u \sim C_5$), by Brook's Theorem $\chi(G-u) \leq \Delta(G-u) < 7 \leq \Delta$. Thus always $\chi(G-u) \leq \Delta-1$ and in fact, $\chi(G-u) = \Delta-1$ and deg $v \geq \Delta-1$ $\forall v \in V(G)$.

Let $Q \subseteq V(G)$ be s.t. $<Q>$ is a maximum clique in G. Let $u \in Q$ be s.t. deg $u = \max_{v \in Q}$ deg $v$. Let $S = \{1,..., \Delta\}$ be a $\Delta$-coloring of G s.t. u is colored $\Delta$ and vertices in Q are colored $1,..., |Q|-1$. Every vertex v of N(u) with a unique color say i has at least one j-vertex $j \neq i$ (else color v by j and u by i). Also as G is $3K_1$-free, $<V(G)-\overline{N(u)}>$ is complete and hence by maximality $|Q| \geq |V(G)-\overline{N(u)}|$     **I**

**Case 1:** deg u = $\Delta-1$ for every maximal clique $<Q>$ in G and a vertex u of maximum degree in Q.
Then deg $v = \Delta-1$ $\forall v \in Q$. As G is $3K_1$-free and deg $u \geq 6$, $|Q| \geq 4$. Let A, C, D $\in Q$ be colored 1, 2, 3 resply. Clearly every j-vertex v of Q has a unique i-vertex in N(v) where $i \neq j$ (else N(v) has color i missing in $\overline{N(v)}$. Color v by i, u by j). As $\Delta > \omega$, |N(u)-Q| ≥ 1. W.l.g. let B ∈ N(u)-Q be colored ω and A ∈ Q be s.t. AB ∉ E(G). Let R be a component containing A s.t. vertices in R are colored either 1 or ω. Then B ∈R (else alter colors in R and color u by 1). Let P = {A, E, F, B} be a 1-ω path in R. F is the only 1-vertex of B and FV ∉ E(G) for any V ∈ Q. Then F has a r-vertex adjacent in V(G)-$\overline{N(u)}$ ∀ r, 2 ≤ r ≤ |Q|-1 (else color F by r, B by 1 and u by ω) . By **I**, E is adjacent to all these |Q|-2 vertices and |V(G)-$\overline{N(u)}$| ≥ |Q|. Hence by **I**, |V(G)-$\overline{N(u)}$| = |Q|, |N(u)-Q| = 1 (else if N(u)-Q has one more vertex, then it is non-adjacent to some vertex in Q and |V(G)-$\overline{N(u)}$| > |Q|) and B is non-adjacent to all vertices of Q ⇒ E is adjacent to all vertices of Q and by assumption Δ-1= deg E ≥ 2(|Q|-1) = 2(Δ-1) ⇒ Δ = 1, a contradiction.

**Case 2:** Let deg u = Δ.
Let A, B, C ∈ Q be colored 1, 2, 3
First let N(u)-Q have a repeat color and D, E ∈ N(u)-Q be colored 4 and F, G ∈ N(u)-Q be colored 5, 6 resply. W.l.g. let DB ∉ E(G) ⇒ BE ∈ E(G) and w.l.g. let EA ∉ E(G) ⇒ AD ∈ E(G). Now D has a 2-vertex in V(G)-$\overline{N(u)}$ (else color D by 2, A by 4, u by 1). Similarly E has a 1-vertex in V(G)-$\overline{N(u)}$. Also as <Q> is a maximum clique in G, F (G) is non-adjacent to some vertex in Q and V(G)-$\overline{N(u)}$ has a 5-vertex and a 6-vertex. Then by **I**, |V(G)-$\overline{N(u)}$| = |Q| and V(G)-$\overline{N(u)}$ has no 3-vertex. Clearly N(C) = N(u) (else G is (Δ-1)-colorable). But then as G is $3K_1$-free, ω<A, B, D, E, F, G> ≥ 3 and ω(G) ≥ 5, a contradiction.
Next w.l.g. let D, E, F, G ∈ N(u)-Q be colored 1, 4, 5, 6 resply. Now ω<E, F, G> ≥ 2. W.l.g. let EF ∈ E(G). Again each of E, F, G is non-adjacent to B or C (else if say EB, EC ∈ E(G), then by replacing Q with Q-A+E we get the earlier case). Thus |V(G)-$\overline{N(u)}$| ≥ 5, a contradiction.

This proves the result.

**Main Result 2:** If G is $3K_1$-free and Δ ≥ 8, then χ ≤ max{ω, Δ-1}.

Proof: As before let G be a smallest $3K_1$-free graph with Δ ≥ 8 and χ > max{ω, Δ-1}. As before we have χ(G-u) = Δ-1 and deg u ≥ Δ-1 ∀ v ∈ V(G). Let |Q| ≥ 5.

Let Q ⊆ V(G) be s.t. <Q> is a maximum clique in G. Let u ∈ Q be s.t. deg u = $\max\limits_{v \in Q}$ deg v. Let S = {1,..., Δ} be a Δ-coloring of G s.t. u is colored Δ and vertices in Q are colored 1,..., |Q|-1. Every vertex v of N(u) with a unique color say i has at least one j-vertex j ≠ i (else color v by j and u by i). Also as G is $3K_1$-free, <V(G)-$\overline{N(u)}$> is complete and hence by maximality |Q| ≥ |V(G)-$\overline{N(u)}$|      **I**

**Case 1:** deg u = Δ-1 for every maximal clique <Q> in G and a vertex u of maximum degree in Q.
Same proof as in the **Main Result 1** holds good.

**Case 2:** deg u = Δ ≥ 8

**Case 2.1:** N(u)-Q has a repeat color
Let A, B ∈ N(u)-Q be colored 5. W.l.g. let ∃ D colored 2 in Q, s.t. AD ∉ E(G). Then DB ∈ E(G). Also let ∃ C colored 1 in Q, s.t. BC ∉ E(G). Then CA ∈ E(G). Let E, F ∈ Q be colored 3, 4 resply.

**Case 2.1.1:** ∃ a vertex X ∈ Q s.t. XA, XB ∈ E(G). W.l.g. let X = E.

**Case 2.1.1.1:** ∃ a vertex Y ∈ Q-C s.t. YA ∈ E(G) and YB ∉ E(G). W.l.g. let Y = F. As D has at the most one repeat color in $\overline{N(D)}$, w.l.g. let C be the only 1-vertex of D. Now let G be the 2-vertex of A in V(G)- $\overline{N(u)}$ (else color A by 2, C by 5, u by 1) and similarly H, J be the 1-vertex, 4-vertex of B in V(G)- $\overline{N(u)}$ resply. Now AH ∈ E(G) (else color A by 1, C by 5, D by 1, u by 2) ⇒ F is the only 4-vertex of A ⇒ DJ ∈ E(G) (else color A by 4, F by 5, D by 4, u by 2). Thus E is the only 3-vertex of A and D. Also as E has two 5-vetices D is its only 2-vertex. Color E by 2, D by 3, A by 3, C by 5, u by 1, a contradiction.

**Case 2.1.1.2:** ∀ V ∈ Q-{C, D}, VA, VB ∈ E(G).
As before let G be the 2-vertex of A and H be the 1-vertex of B. As E, F have two 5-vertices each, C, D are the only 1, 2 vertices of E and F. Also w.l.g. let E be the only 3-vertex of A ⇒ DJ ∈ E(G) where J is the 3-vertex in V(G)- $\overline{N(u)}$ (else color E by 1, C by 5, A by 3, D by 3, u by 2). Again AK ∈ E(G) where K is the 4-vertex in V(G)- $\overline{N(u)}$ (else color F by 1, C by 5, A by 4, D by 4, u by 2). Thus C (A) is the only 1 (5) vertex of A (C) and C is the only 1-vertex of D. Color C by 5, A by 1, D by 1, u by 2, a contradiction.

**Case 2.1.2:** Every V ∈ Q is adjacent to only one of A or B.

**Case 2.1.2.1:** ∃ two vertices in Q-C say E, F adjacent to say A.
W.l.g. let C, E be the unique 1, 3 vertices of A and E be the unique 3-vertex of D. Color E by 5, A by 3, D by 3, u by 2, a contradiction.

**Case 2.1.2.2:** No two vertices in Q-C (Q-D) are adjacent to A (B).
Then |Q| = 5 (else we get **Case 2.1.2.1**). W.l.g. let EA, FB ∈ E(G) and C, D be the only 1, 2 vertex of A, B resply ⇒ C (D) has another 2 (1) vertex in V(G)- $\overline{N(u)}$ and A (B) has another 3 (4) vertex in V(G)- $\overline{N(u)}$ and |V(G)- $\overline{N(u)}$| ≥ 4. Also as Δ ≥ 8, and |Q| = 5, |N(u)-Q-{A, B}| ≥ 2. Clearly |V(G)- $\overline{N(u)}$| ≥ 6, contrary to **I**.

**Case 2.2:** A ∈ Q and B ∈ N(u)-Q have a common color 1.
Now every vertex say V ∈ R = N(u)-Q-B is non-adjacent to some vertex of Q-u-A (else replace Q by Q-A+V to get **Case 2.1**). Let C, D, E ∈ Q have colors 2, 3, 4 resply.

**Claim:** |Q| ≥ 6
If |Q| = 5, then as Δ ≥ 8, |R| ≥ 4 and ω<R> ≥ 2. W.l.g. let FG ∈ E(G) and FC, GD ∉ E(G). Let F, G, H ∈ R be colored 5, 6, 7 resply. ⇒ V(G)- $\overline{N(u)}$ has a 2-vertex, 3-vertex, 5-vertex, 6-vertex and a 7-vertex. Hence by **I** E is the unique 4-vertex in V(G). As E has at the most one repeat color in N(E), w.l.g. let F be the only 5-vertex of E. Color F by 4, E by 5, C by 4, u by 2, a contradiction.

This proves the **Claim.**

Let L ∈ Q be colored 6 and F ∈ R be colored 5.

**Case 2.2.1:** ∃ a vertex in R non-adjacent to two vertices in Q-A.
W.l.g. let FC, FD ∉ E(G). Then ∃ a 2-5 path P = {C, G, H, F}. Similarly let R = {D, G, J, F} 3-5 path. As G has two 2-vertices and 3-vertices, G has a unique 4-vertex (else some color r is missing in $\overline{N(G)}$, color G by r, C by 5, u by 2). Then GE ∉ E(G) (else color E by 2, C by 5, G by 4, u by 4), F is the only 5-vertex of E and let K ∈ V(G)- $\overline{N(u)}$ be the unique 4-vertex of G.

If |Q| ≥ 7, then let M be the 7-vertex in Q. Now as before G has a unique 6-vertex, 7-vertex and LG, MG ∉ E(G). Also as F has at the most one repeat color in $\overline{N(F)}$ w.l.g. let E, L be the only 4-vertex, 6-vertex of F. Again C has either a unique 4-vertex or 6-vertex. W.l.g. E be the only 4-vertex of C. Color E by 5, F by 4, C by 4, u by 2, a contradiction.

Hence <u>let |Q| = 6</u>. W.l.g. let E be the only 4-vertex of F. If E is the only 4-vertex of C (D), then color E by 5, F by 4, C (D) by 4, u by 2 (3), a contradiction. Hence CK, DK ∈ E(G) where K is the 4-vertex in V(G)-$\overline{N(u)}$. Then L is the only 6-vertex of C, D. If L is the only 6-vertex of F, then color L by 5, F by 6, C by 6, u by 2, a contradiction. Hence F has another 6-vertex M in V(G)-$\overline{N(u)}$. As Δ ≥ 8, and |Q| = 6, ∃ N, P ∈ N(u)-Q colored 7, 8 resply. By **I**, either N or P is the unique 7-vertex or 8-vertex in V(G). W.l.g. let N be the unique 7-vertex in V(G). As N has a unique 2-vertex or 5-vertex, w.l.g. let C be its unique 2-vertex. Color N by 2, C by 7, F by 7, u by 5, a contradiction.

**Case 2.2.2:** Every vertex of R is non-adjacent to exactly one vertex in Q-A.

Let FC ∉ E(G). As F has at most one color repeated in $\overline{N(F)}$ w.l.g. let D, E be the unique 3-vertex and 4-vertex of F. Again w.l.g. let D be the unique 3-vertex of C. Now D has either a unique 2-vertex or 5-vertex. W.l.g. let C be the unique 2-vertex of D. Color D by 2, C by 3, F by 3, u by 5, a contradiction.

This proves the theorem.

**Corollary:** Borodin and Kostochka conjecture is true for $3K_1$-free graphs.